# THE PROBABILITY OF EXCEEDING A HIGH BOUNDARY ON A RANDOM TIME INTERVAL FOR A HEAVY-TAILED RANDOM WALK

By Serguei Foss[1], Zbigniew Palmowski[2] and Stan Zachary

*Heriot–Watt University, University of Wrocław and Utrecht University, and Heriot–Watt University*

We study the asymptotic probability that a random walk with heavy-tailed increments crosses a high boundary on a random time interval. We use new techniques to extend results of Asmussen [*Ann. Appl. Probab.* **8** (1998) 354–374] to completely general stopping times, uniformity of convergence over all stopping times and a wide class of nonlinear boundaries. We also give some examples and counterexamples.

**1. Introduction and main results.** The analysis of random walks with heavy-tailed increments is central to the understanding of many problems in insurance, finance, queueing networks and storage theory. In particular, we are often interested in determining the probability of overcrossing a deterministic curve $\{x + g(n)\}_{n\geq 0}$ as $x$ is allowed to become large.

Thus, in this paper, we consider a sequence $\{\xi_n\}_{n\geq 1}$ of independent identically distributed random variables with distribution function $F$. We assume throughout that $F$ belongs to the class $\mathcal{L}$ of *long-tailed* distribution functions, where a distribution function $G \in \mathcal{L}$ if and only if

(1) $\quad \bar{G}(x) > 0 \quad$ for all $x$, $\quad \lim_{x\to\infty} \frac{\bar{G}(x-h)}{\bar{G}(x)} = 1 \quad$ for all fixed $h > 0$.

Here $\bar{G}$ denotes the tail distribution given by $\bar{G}(x) = 1 - G(x)$. We further assume throughout that the distribution $F$ has a finite mean $m_F = \mathbf{E}\xi_1$. Without loss of generality (see below), we assume

$$m_F = 0.$$

Received May 2004; revised October 2004.
[1]Supported in part by the EPSRC Grant No GR/R58765/01.
[2]Supported in part by the NWO Grant No 613.000.310 and Heriot–Watt University.
*AMS 2000 subject classifications.* Primary 60J10; secondary 60K25.
*Key words and phrases.* Ruin probability, subexponential distributions, random walk, boundary.

Define the random walk $\{S_n\}_{n\geq 0}$ by

$$S_0 = 0, \qquad S_n = \sum_{i=1}^{n} \xi_i, \qquad n \geq 1.$$

For any nonnegative function $g$ on $\mathbb{Z}_+$, define also the process $\{S_n^g\}_{n\geq 0}$ by

$$S_n^g = S_n - g(n), \qquad n \geq 0.$$

The process $\{S_n^g\}_{n\geq 0}$ is investigated in *nonlinear* renewal theory (see [22]), and also in many other examples in probability and queueing theory (see, e.g., [1, 7, 20, 21, 24]). Note also that any subadditive functional of a random walk is of this form—see [8].

For $n \geq 0$, let

$$M_n^g = \max_{0 \leq i \leq n} S_i^g.$$

Similarly, for any stopping time $\sigma$ for the random walk $\{S_n\}_{n\geq 0}$ (i.e., for any random variable $\sigma$ taking values in $\mathbb{Z}_+ \cup \{\infty\}$ such that, for all $n \geq 0$, the event $\{\sigma \leq n\}$ is independent of $\{\xi_{n+1}, \xi_{n+2}, \ldots\}$), let

$$M_\sigma^g = \max_{0 \leq i \leq \sigma} S_i^g.$$

Define also the decreasing function $H_\sigma^g$ by

$$H_\sigma^g(x) = \sum_{n \geq 1} \mathbf{P}(\sigma \geq n) \bar{F}(x + g(n)).$$

Note that the function $H_\sigma^g$ is monotone decreasing in $g$ [i.e., if $g_1(n) \geq g_2(n)$ for all $n$, then $H_\sigma^{g_1}(x) \leq H_\sigma^{g_2}(x)$ for all $x$] and monotone increasing in $\sigma$ [i.e., if $\sigma_1 \geq \sigma_2$ a.s., then $H_{\sigma_1}^g(x) \geq H_{\sigma_2}^g(x)$ for all $x$]. Note also that, since $F$ has a finite mean, $H_\sigma^g$ is finite for all $\sigma$ and all $g$ such that $g(n) \geq cn$ for some $c > 0$; further, since $F \in \mathcal{L}$, an elementary truncation argument along the lines of the proof of Lemma 1(i) shows that, for any $\sigma$ such that $\mathbf{E}\sigma < \infty$ and nonnegative function $g$, $H_\sigma^g(x)$ is finite for all $x$ and

(2) $$H_\sigma^g(x) = (1+o(1))\mathbf{E}\sigma \bar{F}(x) \qquad \text{as } x \to \infty.$$

We are interested in the asymptotic distribution of $M_\sigma^g$ for a general stopping time $\sigma$ (which need not be a.s. finite). In particular, we are interested in obtaining conditions under which

(3) $$\mathbf{P}(M_\sigma^g > x) \geq (1+o(1))H_\sigma^g(x) \qquad \text{as } x \to \infty,$$

and in obtaining (stronger) conditions under which

(4) $$\mathbf{P}(M_\sigma^g > x) = (1+o(1))H_\sigma^g(x) \qquad \text{as } x \to \infty,$$

in each case with uniformity over suitable classes of stopping times $\sigma$ and functions $g$. [We shall say, e.g., that the result (3) holds with uniformity over all $\sigma$ and all $g$—in appropriate classes—if and only if there exists a function $\delta$ on $\mathbb{R}_+$ such that $\delta(x) \to 0$ as $x \to \infty$ and $\mathbf{P}(M_\sigma^g > x) \geq (1-\delta(x))H_\sigma^g(x)$ for all $x \in \mathbb{R}_+$ and for all $\sigma$ and all $g$.]

The event $\{M_\sigma^g > x\}$ may be reinterpreted as the event that the random walk $\{S_n\}_{n\geq 0}$ crosses the (arbitrary) increasing boundary $\{x+g(n)\}_{n\geq 0}$ by the stopping time $\sigma$. The intuitive interpretation of the relation (4), in particular, is that, for $x$ very large, the only significant way in which the random walk can cross this boundary is that it remains close to its mean zero up to some time $n$ when, with probability $\bar{F}(x+g(n))$, it jumps above $x + g(n)$. This property is the "principle of one big jump" and is characteristic of the subexponential property (see below) which we shall in general require (at a minimum) of $F$ in order to obtain conditions for (4) to hold.

Our results below are also applicable to random walks whose increments have a nonzero mean: it is clearly sufficient to make the obvious shift transformation. In particular, by considering, for $c > 0$, the function $g(n) = cn$, the results include as a special case those for the maximum on a random interval of a random walk with drift $-c$. The results obtained in this case both generalize and extend earlier results of Asmussen [2] and Foss and Zachary [14]. We give a more detailed discussion of this below.

In order to state our results, we require some further definitions. A distribution function $G$ on $\mathbb{R}_+$ is *subexponential* if and only if $\bar{G}(x) > 0$ for all $x$ and

(5) $$\lim_{x\to\infty} \overline{G^{*2}}(x)/\bar{G}(x) = 2$$

(where $G^{*2}$ is the convolution of $G$ with itself). More generally, a distribution function $G$ on $\mathbb{R}$ is subexponential if and only if $G^+$ is subexponential, where $G^+ = G\mathbf{I}_{\mathbb{R}_+}$ and $\mathbf{I}_{\mathbb{R}_+}$ is the indicator function of $\mathbb{R}_+$. It is known that the subexponentiality of a distribution depends only on its (right) tail, and that a subexponential distribution is long-tailed. We let $\mathcal{S}$ denote the class of subexponential distributions, so that, in particular, $\mathcal{S} \subset \mathcal{L}$.

A distribution function $G$ on $\mathbb{R}$ belongs to the class $\mathcal{S}^*$ introduced by [16] if and only if $\bar{G}(x) > 0$ for all $x$ and

(6) $$\int_0^x \bar{G}(x-y)\bar{G}(y)\,dy \sim 2m_{G^+}\bar{G}(x) \qquad \text{as } x \to \infty,$$

where

$$m_{G^+} = \int_0^\infty \bar{G}(x)\,dx$$

is the mean of $G^+$. It is again known that the property $G \in \mathcal{S}^*$ depends only on the tail of $G$. Further, if $G \in \mathcal{S}^*$ then $G \in \mathcal{S}$, and also $G^s \in \mathcal{S}$, where

$$\overline{G^s}(x) = \min\left(1, \int_x^\infty \bar{G}(t)\,dt\right)$$

is the integrated, or *second-tail*, distribution function determined by $G$—see [16].

Let $\mathcal{T}$ be the class of all stopping times for the random walk $\{S_n\}_{n\geq 0}$. For any stopping time $\varphi$, let

$$\mathcal{T}_\varphi = \{\sigma \in \mathcal{T} : \sigma \leq \varphi \text{ a.s.}\}.$$

In particular, for any integer $N > 0$, $\mathcal{T}_N$ is the class of stopping times almost surely bounded by $N$.

For any constant $c$ (we shall primarily be interested in $c \geq 0$), let $\mathcal{G}_c$ be the class of nonnegative functions $g$ satisfying

(7) $\qquad g(1) \geq c, \qquad g(n+1) \geq g(n) + c, \qquad n \geq 1.$

In particular, $\mathcal{G}_0$ is the class of nonnegative nondecreasing functions on $\mathbb{Z}_+$. Note also that the class $\mathcal{G}_c$ is monotone decreasing in $c$.

As a preliminary result, we prove the following theorem, which relates to bounded stopping times.

THEOREM 1. (i) *Suppose that $F \in \mathcal{L}$. Then, given any integer $N > 0$, the result* (3) *holds uniformly over all $\sigma \in \mathcal{T}_N$ and all $g \in \mathcal{G}_0$.*

(ii) *Suppose, additionally, that $F \in \mathcal{S}$. Then, given any integer $N > 0$, the result* (4) *holds uniformly over all $\sigma \in \mathcal{T}_N$ and all $g \in \mathcal{G}_0$.*

Our main result is then Theorem 2.

THEOREM 2. (i) *Suppose that $F \in \mathcal{L}$. Then, given any $c > 0$, the result* (3) *holds uniformly over all $\sigma \in \mathcal{T}$ and all $g \in \mathcal{G}_c$.*

(ii) *Suppose, additionally, that $F \in \mathcal{S}^*$. Then, given any $c > 0$, the result* (4) *holds uniformly over all $\sigma \in \mathcal{T}$ and all $g \in \mathcal{G}_c$.*

We have stated these results under those conditions which appear to us most natural. There are some obvious extensions which are immediate from the condition $F \in \mathcal{L}$ which we assume throughout. This condition implies that also $H_\sigma^g \in \mathcal{L}$ with uniformity of convergence in the definition (1) over all stopping times $\sigma$ and nonnegative functions $g$. Thus, for example, for any $c$ for which one of the results of Theorems 1 or 2 holds, and for any fixed $d > 0$, we may expand the corresponding class $\mathcal{G}_c$ to include any function $g$ such that $g' \leq g \leq g' + d$ for some function $g' \in \mathcal{G}_c$—since then $H_\sigma^{g'}(x)/H_\sigma^{g'+d}(x) \to 1$ as $x \to \infty$ with the required uniformity properties. One consequence of this observation is that we may, in either of the results of Theorem 1, replace $\mathcal{G}_0$ by $\mathcal{G}_c$ for any $c \in \mathbb{R}$. That we may not, in general, even for a single bounded stopping time $\sigma$, obtain the results of Theorem 1 with uniformity over all nonnegative functions $g$ is shown by Example 1 of Section 3. See also that section for further discussion and comments.

We now discuss briefly our main result, which is part (ii) of Theorem 2. Consider first the slightly weaker condition $F^s \in \mathcal{S}$, and the case where the function $g$ is given by $g(n) = cn$ for some $c > 0$, and the stopping time $\sigma = \infty$. The conclusion (4) is then equivalent to the well-known result of Veraverbeke [19] for the asymptotic distribution of the maximum of the random walk $\{S_n - cn\}_{n \geq 0}$ with drift $-c$ (see the Appendix for conditions under which the function $H_\sigma^g$ has a tail equivalent integral representation). See also [12] and [11]. Now assume that $F \in \mathcal{S}^*$. In the case where the function $g$ is again given by $g(n) = cn$ for $c > 0$, and the stopping time $\sigma$ has a finite mean, it follows from (2) that the conclusion (4) is equivalent to

$$(8) \qquad \mathbf{P}(M_\sigma^g > x) = (1 + o(1))\mathbf{E}\sigma \bar{F}(x) \qquad \text{as } x \to \infty.$$

(Again the event $\{M_\sigma^g > x\}$ is most naturally interpreted in relation to the random walk $\{S_n - cn\}_{n \geq 0}$.) Asmussen [2] proved the result (8) for the stopping time $\sigma = \tau_c \equiv \min\{n \geq 1 : S_n - cn \leq 0\}$ (see also [15]). Foss and Zachary [14] extended the result (8) to a general stopping time $\sigma$, and showed also the necessity (for a general stopping time) of the condition $F \in \mathcal{S}^*$. However, in the case $\mathbf{E}\sigma = \infty$ (which occurs naturally in many applications—see, e.g., Example 3 of Section 3), the result (8) simply asserts that $\mathbf{P}(M_\sigma^g > x)/\bar{F}(x) \to \infty$ as $x \to \infty$ and does not give the asymptotic form of the tail of the distribution of $M_\sigma^g$. Nor, as may be deduced from the results of the present paper, does the result (8) hold with uniformity even over all finite stopping times $\sigma$. In the present paper we obtain the correct asymptotics in the case $\mathbf{E}\sigma = \infty$, we extend our results to arbitrary boundaries $g$, and we give these results in such a form that in each case we obtain uniformity of convergence over all stopping times $\sigma$ (which need not be a.s. finite) and over suitably wide classes of functions $g$. This uniformity corresponds to the naturalness of the condition (4), as discussed above, and of course guarantees the quality of the asymptotic results, notably over all $\sigma$. In particular, Theorem 2 thus unifies the earlier, quite distinct, results for the cases $\mathbf{E}\sigma < \infty$ and $\sigma = \infty$ a.s. (with $g$ linear in each case).

Further, in the present paper we take Asmussen's result (8) for $g(n) = cn$ and $\sigma = \tau_c$ as a starting point and use new and direct arguments to obtain our results for general stopping times $\sigma$ and classes of functions $g$. (Notably, we make no further use, beyond its requirement for Asmussen's result, of the condition $F \in \mathcal{S}^*$.) Denisov [9] has recently given a very simple proof of (8) for $g(n) = cn$ and $\sigma = \tau_c$. This, taken with the present paper, now yields a relatively simple and direct treatment of all our results.

We note also here that, in the case where the stopping time $\sigma$ is *independent* of $\{S_n\}_{n \geq 0}$ and the function $g$ is given by $g(n) = cn$ for $c > 0$, that the result (4) holds with uniformity over all such $\sigma$ follows from the results of Korshunov [17]—see the comments on this in [10].

In Section 2 we prove our main results, giving parallel developments of the lower and upper bounds so as to identify carefully the conditions required for each. We prove our results successively for bounded stopping times (Theorem 1 above), stopping times bounded by a stopping time with a finite mean (for the upper bound we require the stopping time $\tau_c$ identified above) and for quite general stopping times (Theorem 2 above).

In Section 3 we give various examples to show the applicability of the results, together with counterexamples to show what goes wrong when we drop the conditions of our theorems.

The Appendix gives a simple integral representation, under appropriate conditions, of the function $H_\sigma^g$.

## 2. Proofs.

PROOF OF THEOREM 1. Since the results are trivial in the case $\sigma = 0$ a.s. and by otherwise conditioning on the event $\{\sigma > 0\}$, we may assume throughout without loss of generality that $\sigma \geq 1$ a.s.

Since $F \in \mathcal{L}$ throughout, we may choose a function $h : \mathbb{R}_+ \to \mathbb{R}_+$ such that

$$(9) \quad h(x) \leq x \quad \text{for all } x \geq 0,$$

$$(10) \quad h \text{ is increasing}, \quad h(x) \to \infty \text{ as } x \to \infty,$$

$$(11) \quad \frac{\bar{F}(x - h(x))}{\bar{F}(x)} \to 1 \quad \text{as } x \to \infty.$$

(This follows from the condition $F \in \mathcal{L}$ by allowing the function $h$ to increase sufficiently slowly—see [14].)

Note that the results of both parts of the theorem are trivial in the case $N = 1$. Given any integer $N \geq 2$, consider any stopping time $\sigma \in \mathcal{T}_N$ and any function $g \in \mathcal{G}_0$. Then, for $x \geq 0$,

$$
\begin{aligned}
\mathbf{P}(M_\sigma^g > x) &= \sum_{n=1}^{N} \mathbf{P}(\sigma \geq n, M_{n-1}^g \leq x, S_n^g > x) \\
&= \sum_{n=1}^{N} \mathbf{P}(\sigma \geq n, M_{n-2}^g \leq x, \\
&\qquad\qquad S_{n-1} < -h(x + g(n-1)), S_n^g > x) \\
&\quad + \sum_{n=1}^{N} \mathbf{P}(\sigma \geq n, M_{n-2}^g \leq x,
\end{aligned}
$$

(12)

$$S_{n-1} \in [-h(x+g(n-1)), h(x+g(n-1))],$$
$$S_n^g > x)$$

$$+ \sum_{n=1}^{N} \mathbf{P}(\sigma \geq n, M_{n-2}^g \leq x,$$
$$S_{n-1} \in (h(x+g(n-1)), x+g(n-1)), S_n^g > x),$$

where, for $n=1$, we may take $M_{n-2}^g = 0$.

Since, for $g \in \mathcal{G}_0$,

$$\mathbf{P}(M_{n-2}^g > x) \leq \mathbf{P}(M_{n-2} > x) \to 0,$$

as $x \to \infty$, and, from (10),

$$\mathbf{P}(S_{n-1} \notin [-h(x+g(n-1)), h(x+g(n-1))]) \leq \mathbf{P}(S_{n-1} \notin [-h(x), h(x)]) \to 0$$

as $x \to \infty$, it follows that, for $1 \leq n \leq N$,

(13)
$$\mathbf{P}(\sigma \geq n, M_{n-2}^g \leq x, S_{n-1} \in [-h(x+g(n-1)), h(x+g(n-1))])$$
$$= \mathbf{P}(\sigma \geq n) + o(1)$$

as $x \to \infty$, uniformly over all $\sigma \in \mathcal{T}_N$ and $g \in \mathcal{G}_0$. Further, it follows from (11) that, for any $n$,

(14)
$$\frac{\bar{F}(x+g(n) \pm h(x+g(n)))}{\bar{F}(x+g(n))} \to 1 \qquad \text{as } x \to \infty,$$

uniformly over all $g \in \mathcal{G}_0$. Since also, for any $n$, $h(x+g(n-1)) \leq h(x+g(n))$, it follows from (13) and (14) that

(15)
$$\sum_{n=1}^{N} \mathbf{P}(\sigma \geq n, M_{n-2}^g \leq x,$$
$$S_{n-1} \in [-h(x+g(n-1)), h(x+g(n-1))], S_n^g > x)$$
$$= (1+o(1)) \sum_{n=1}^{N} (\mathbf{P}(\sigma \geq n) + o(1)) \bar{F}(x+g(n))$$
$$= (1+o(1)) H_\sigma^g(x) + o(\bar{F}(x+g(1)))$$
$$= (1+o(1)) H_\sigma^g(x)$$

as $x \to \infty$, uniformly over all $\sigma \in \mathcal{T}_N$ and $g \in \mathcal{G}_0$, where the final line in (15) follows since $\sigma \geq 1$ a.s. Since the first and third terms on the right-hand side of (12) are positive, the result (i) of the theorem now follows from (12) and (15).

To prove (ii), we suppose that $F \in \mathcal{S}$. We require to show that (4) holds uniformly over all $\sigma \in \mathcal{T}_N$ and $g \in \mathcal{G}_0$. From (12) and (15), it is sufficient to show that the first and third terms on the right-hand side of (12) are each $o(H_\sigma^g(x))$ as $x \to \infty$, again uniformly over all $\sigma \in \mathcal{T}_N$ and $g \in \mathcal{G}_0$. That this is true for the first of these terms follows since, for each $n$,

$$\mathbf{P}(\sigma \geq n, M_{n-2}^g \leq x, S_{n-1} < -h(x+g(n-1)), S_n^g > x)$$
$$\leq \mathbf{P}(S_{n-1} < -h(x+g(n-1)))\bar{F}(x+g(n))$$
$$\leq \mathbf{P}(S_{n-1} < -h(x))\bar{F}(x+g(1))$$
$$= o(H_\sigma^g(x))$$

as $x \to \infty$ [from (10) and since $\sigma \geq 1$ a.s.] with the required uniformity.

In the case where $\sigma$ is identically equal to $N$ and $g$ is identically equal to 0, it is a standard result that

$$\mathbf{P}(M_\sigma^g > x) = \mathbf{P}\left(\max_{0 \leq n \leq N} S_n > x\right)$$
$$= (1+o(1))N\bar{F}(x) \qquad \text{as } x \to \infty,$$
(16)

(see [11]). Since in this case $H_\sigma^g(x) = N\bar{F}(x)$, it follows from (12), (15) and (16) that, for $1 \leq n \leq N$,

(17) $\mathbf{P}(S_i \leq x, i \leq n-2; S_{n-1} \in (h(x), x]; S_n > x) = o(\bar{F}(x))$ as $x \to \infty$.

For general $\sigma \in \mathcal{T}_N$ with $\sigma \geq 1$ a.s. and $g \in \mathcal{G}_0$, it follows since $g$ is nondecreasing that the third term on the right-hand side of (12) is bounded above by

$$\sum_{n=1}^{N} \mathbf{P}(S_i \leq x + g(n-1), i \leq n-2,$$
$$S_{n-1} \in (h(x+g(n-1)), x+g(n-1)], S_n > x+g(n-1)).$$

From (17), the $n$th term in the above sum is $o(\bar{F}(x+g(n-1)))$, and so also (since $\sigma \geq 1$ a.s. and $g \in \mathcal{G}_0$) the sum is $o(H_\sigma^g(x))$, as $x \to \infty$, uniformly over all such $\sigma$ and $g$ as required. □

REMARK 1. In Section 3 we give examples which show that we may not, in general, drop the condition that $g$ be nondecreasing.

The proof of our main result, Theorem 2, requires the separate derivation of upper and lower bounds for $\mathbf{P}(M_\sigma^g > x)$. In Lemma 1 below, we first establish these bounds for classes of stopping times intermediate between those of Theorems 1 and 2.

For any $a > 0$, define the stopping time $\tau_a = \min\{n \geq 1 : S_n < an\}$. Note that, since $F$ has mean 0, $\mathbf{E}\tau_a$ is finite. For any $a > 0$, define also the function $\bar{a}$ on $\mathbb{Z}_+$ by $\bar{a}(n) = an$.

LEMMA 1. (i) *Given any stopping time $\varphi$ such that $\mathbf{E}\varphi < \infty$, the result* (3) *holds uniformly over all $\sigma \in \mathcal{T}_\varphi$ and all $g \in \mathcal{G}_0$.*

(ii) *Suppose that $F \in \mathcal{S}^*$. Then, given any $c > 0$, the result* (4) *holds uniformly over all $\sigma \in \mathcal{T}_{\tau_c}$ and all $g \in \mathcal{G}_c$.*

PROOF. In the proofs of both (i) and (ii), we may again assume without loss of generality, as in the proof of Theorem 1, that $\sigma \geq 1$ a.s. Thus, given $\varphi$ such that $\mathbf{E}\varphi < \infty$, for any $\sigma \in \mathcal{T}_\varphi$ with $\sigma \geq 1$ a.s. and $g \in \mathcal{G}_0$, and for any integer $N > 0$ and all $x > 0$,

$$H_\sigma^g(x) - H_{\sigma \wedge N}^g(x) = \sum_{n > N} \mathbf{P}(\sigma \geq n) \bar{F}(x + g(n))$$

$$\leq \bar{F}(x + g(1)) \sum_{n > N} \mathbf{P}(\sigma \geq n)$$

$$\leq H_\sigma^g(x) \sum_{n > N} \mathbf{P}(\sigma \geq n)$$

$$\leq H_\sigma^g(x) \sum_{n > N} \mathbf{P}(\varphi \geq n).$$

Hence, using Theorem 1(i) applied to the stopping time $\sigma \wedge N$, there exists a function $\varepsilon_N$, which is independent of $\sigma$ and $g$, such that $\varepsilon_N(x) \to 0$ as $x \to \infty$ and, for $\sigma$ and $g$ as above and for $x > 0$,

$$\mathbf{P}(M_\sigma^g > x) \geq \mathbf{P}(M_{\sigma \wedge N}^g > x)$$

$$\geq (1 - \varepsilon_N(x)) H_{\sigma \wedge N}^g(x)$$

$$\geq (1 - \varepsilon_N(x)) H_\sigma^g(x) \left(1 - \sum_{n > N} \mathbf{P}(\varphi \geq n)\right).$$

Since $\mathbf{E}\varphi < \infty$, it now follows that

(18) $$\mathbf{P}(M_\sigma^g > x) \geq (1 - \varepsilon'_N(x)) H_\sigma^g(x)$$

for some positive function $\varepsilon'_N$, again independent of $\sigma$ and $g$, such that

$$\lim_{N \to \infty} \lim_{x \to \infty} \varepsilon'_N(x) = 0.$$

This latter condition implies that (for any such sequence of functions $\{\varepsilon'_N\}_{N \geq 1}$) there exists an integer-valued function $N$ on $\mathbb{R}_+$ such that $\lim_{x \to \infty} \varepsilon'_{N(x)}(x) = 0$. Hence, from (18), we have the required result (3) with the required uniformity over $\sigma \in \mathcal{T}_\varphi$ and $g \in \mathcal{G}_0$.

To prove (ii), we suppose that $F \in \mathcal{S}^*$ and that $c > 0$. Consider first the stopping time $\sigma = \tau_c$ and the function $g = \bar{c}$. For integer $N > 0$, it follows from the result of Asmussen [2] referred to in the Introduction—see also [3], Chapter X, Theorem 9.4—that, as $x \to \infty$,

$$
\begin{aligned}
\mathbf{P}(M^{\bar{c}}_{\tau_c} > x) &= (1+o(1))\mathbf{E}\tau_c \bar{F}(x+c) \\
&= (1+o(1))(\mathbf{E}(\tau_c \wedge N) + \mathbf{E}(\tau_c - N)^+)\bar{F}(x+c) \\
&= (1+o(1))(H^{\bar{c}}_{\tau_c \wedge N}(x) + \mathbf{E}(\tau_c - N)^+ \bar{F}(x+c)),
\end{aligned}
\tag{19}
$$

where (19) follows since $F$ is long-tailed. Since $\mathcal{S}^* \subset \mathcal{S}$, it follows also from Theorem 1(ii) that

$$
\mathbf{P}(M^{\bar{c}}_{\tau_c \wedge N} > x) = (1+o(1))H^{\bar{c}}_{\tau_c \wedge N}(x) \qquad \text{as } x \to \infty. \tag{20}
$$

Since also $H^{\bar{c}}_{\tau_c \wedge N}(x) \leq N\bar{F}(x+c)$, it follows from (19) and (20) that

$$
\begin{aligned}
\mathbf{P}(M^{\bar{c}}_{\tau_c \wedge N} &\leq x, M^{\bar{c}}_{\tau_c} > x) \\
&= \mathbf{P}(M^{\bar{c}}_{\tau_c} > x) - \mathbf{P}(M^{\bar{c}}_{\tau_c \wedge N} > x) \\
&= (1+o(1))\mathbf{E}(\tau_c - N)^+ \bar{F}(x+c) \qquad \text{as } x \to \infty.
\end{aligned}
\tag{21}
$$

We now prove (4) for any $\sigma \in \mathcal{T}_{\tau_c}$ and $g \in \mathcal{G}_c$. For $n \geq 1$, let $d_n = g(n) - cn$. Fix any integer $N > 0$. Then, for $x > 0$,

$$
\mathbf{P}(M^g_{\sigma \wedge N} \leq x, M^g_\sigma > x) \leq \mathbf{P}(M^g_{\tau_c \wedge N} \leq x, M^g_{\tau_c} > x) \tag{22}
$$

$$
\leq \mathbf{P}(M^{\bar{c}}_{\tau_c \wedge N} \leq x + d_N, M^{\bar{c}}_{\tau_c} > x + d_N) \tag{23}
$$

$$
\leq (1+o(1))\mathbf{E}(\tau_c - N)^+ \bar{F}(x+g(1)), \tag{24}
$$

uniformly over all such $\sigma$ and $g$, where (22) follows by consideration of sample paths, while (23) follows since the condition $g \in \mathcal{G}_c$ implies that $d_n$ is nondecreasing in $n$, and finally, (24) follows from (21) on noting that $d_N \geq d_1 = g(1) - c$. Hence, from (24), using Theorem 1(ii) again and noting that, for all $x > 0$, $\bar{F}(x + g(1)) \leq H^g_\sigma(x)$, we have that, as $x \to \infty$,

$$
\mathbf{P}(M^g_\sigma > x) \leq \mathbf{P}(M^g_{\sigma \wedge N} > x) + (1+o(1))\mathbf{E}(\tau_c - N)^+ \bar{F}(x+g(1))
$$

$$
\leq (1 + \mathbf{E}(\tau_c - N)^+ + o(1))H^g_\sigma(x),
$$

uniformly over all $\sigma$ and $g$ as above. Since $\mathbf{E}(\tau_c - N)^+ \to 0$ as $N \to \infty$, we conclude, as in the final part of the proof of part (i) above, that

$$
\mathbf{P}(M^g_\sigma > x) \leq (1+o(1))H^g_\sigma(x) \qquad \text{as } x \to \infty,
$$

again uniformly over all $\sigma$ and $g$ as above. The required result (4) now follows on using also part (i) of the lemma. □

REMARK 2. Note that the result of Asmussen used in the above lemma requires $F \in \mathcal{S}^*$. This is the only point in the argument of the present paper in which this condition is explicitly used.

The proof of the lower bound in Theorem 2 is by consideration of repeated upcrossings by $\{S_n\}_{n\geq 0}$ of boundaries of slope $-a < 0$, while the proof of the upper bound is by consideration of repeated downcrossings of boundaries of slope $a > 0$. In each case $a$ is then allowed to tend to 0. Each argument requires an application of Lemma 1 to the random walk "restarted" at these upcrossing or downcrossing times. We give this in Corollary 1 below, which is stated in a form carefully adapted to its subsequent use.

For any a.s. finite stopping time $\varphi$ and any $a > 0$, define the further stopping time

$$\rho_a^\varphi = \varphi + \min\{n \geq 1 : S_{\varphi+n} - S_\varphi > -an\}.$$

Note that, since $F$ has mean 0, $\rho_a^\varphi$ is a.s. finite.

Similarly, for any a.s. finite stopping time $\varphi$ and $a > 0$, define the further stopping time

$$\tau_a^\varphi = \varphi + \min\{n \geq 1 : S_{\varphi+n} - S_\varphi \leq an\}.$$

Note again that $\tau_a^\varphi$ is a.s. finite.

COROLLARY 1. (i) *Given any $a > 0$, there exists a function $\gamma_a$ on $\mathbb{R}_+$ such that $\lim_{x\to\infty} \gamma_a(x) = 0$ and*

$$\begin{aligned}(25)\quad &\mathbf{P}(\exists n : \varphi < n \leq \sigma \wedge \rho_a^\varphi, S_n^g - S_\varphi^{-\bar{a}} > x) \\ &\qquad \geq (1 - \gamma_a(x)) \sum_{n \geq 1} \mathbf{P}(\varphi < n \leq \sigma \wedge \rho_a^\varphi) \bar{F}(x + g(n) + an),\end{aligned}$$

*for all $x > 0$, all a.s. finite stopping times $\varphi$ and all $\sigma \in \mathcal{T}$ and $g \in \mathcal{G}_0$.*

(ii) *Suppose that $F \in \mathcal{S}^*$. Then, given any $a > 0$, there exists a function $\delta_a$ on $\mathbb{R}_+$ such that $\lim_{x\to\infty} \delta_a(x) = 0$ and*

$$\begin{aligned}(26)\quad &\mathbf{P}(\exists n : \varphi < n \leq \sigma \wedge \tau_a^\varphi, S_n^g - S_\varphi^{\bar{a}} > x) \\ &\qquad \leq (1 + \delta_a(x)) \sum_{n \geq 1} \mathbf{P}(\varphi < n \leq \sigma \wedge \tau_a^\varphi) \bar{F}(x + g(n) - an),\end{aligned}$$

*for all $x > 0$, all a.s. finite stopping times $\varphi$ and all $\sigma \in \mathcal{T}$ and $g \in \mathcal{G}_a$.*

PROOF. We first prove (i). Fix $a > 0$. Note that the stopping time $\rho_a \equiv \rho_a^0 \equiv \min\{n \geq 1 : S_n > -an\}$ has a finite mean. It follows from Lemma 1(i) that there exists a function $\gamma_a$ on $\mathbb{R}_+$ with $\lim_{x\to\infty} \gamma_a(x) = 0$ and such that, for any $\sigma \in \mathcal{T}$ and $g \in \mathcal{G}_0$, and all $x > 0$,

$$\begin{aligned}(27)\quad &\mathbf{P}(\exists n : 0 < n \leq \sigma \wedge \rho_a, S_n^g > x) \\ &\qquad \geq (1 - \gamma_a(x)) \sum_{n \geq 1} \mathbf{P}(n \leq \sigma \wedge \rho_a) \bar{F}(x + g(n)).\end{aligned}$$

Now given $\sigma$ and $g$ as above and any stopping time $\varphi$, to prove (25), we may assume without loss of generality that $\varphi = m$ for some constant $m$ (for otherwise we may condition on each possible value $m$ of $\varphi$, and note that the function $\gamma_a$ is independent of $m$). Thus, consider the random walk $\{S'_n\}_{n \geq 0}$ given by $S'_n = S_{m+n} - S_m$. We have $\rho^m_a - m = \rho'_a$, where $\rho'_a = \min\{n \geq 1 : S'_n > -an\}$, and so the application of (27) to the random walk $\{S'_n\}$, the stopping time $\sigma' = 0 \vee (\sigma - m)$ (for $\{S'_n\}$) and the function $g' \in \mathcal{G}_0$ given by $g'(n) = g(m + n) + am$ gives, for $x > 0$,

$$\mathbf{P}(\exists n : 0 < n \leq (\sigma \wedge \rho^m_a) - m, S^g_{m+n} - S^{-\bar{a}}_m > x)$$
$$= \mathbf{P}(\exists n : 0 < n \leq \sigma' \wedge \rho'_a, S'_n > x + g'(n))$$
(28) $$\geq (1 - \gamma_a(x)) \sum_{n \geq 1} \mathbf{P}(n \leq (\sigma \wedge \rho^m_a) - m) \bar{F}(x + g(m + n) + am)$$
$$\geq (1 - \gamma_a(x)) \sum_{n \geq 1} \mathbf{P}(n \leq (\sigma \wedge \rho^m_a) - m)$$
$$\times \bar{F}(x + g(m + n) + a(m + n)),$$

where the last line follows since $a > 0$. Replace $n$ by $n - m$ in (28) to obtain

$$\mathbf{P}(\exists n : m < n \leq \sigma \wedge \rho^m_a, S^g_n - S^{-\bar{a}}_m > x)$$
$$\geq (1 - \gamma_a(x)) \sum_{n \geq m+1} \mathbf{P}(n \leq \sigma \wedge \rho^m_a) \bar{F}(x + g(n) + an)$$
$$= (1 - \gamma_a(x)) \sum_{n \geq 1} \mathbf{P}(m < n \leq \sigma \wedge \rho^m_a) \bar{F}(x + g(n) + an),$$

which is (25) with $\varphi = m$ as required.

The proof of (ii) is similar to that of (i) with only minor variations. Thus, we suppose that $F \in \mathcal{S}^*$, and fix $a > 0$. It follows from Lemma 1(ii) that there exists a function $\delta_a$ on $\mathbb{R}_+$ with $\lim_{x \to \infty} \delta_a(x) = 0$ and such that, for any $\sigma \in \mathcal{T}$, any $g \in \mathcal{G}_a$, and all $x > 0$,

(29) $$\mathbf{P}(\exists n : 0 < n \leq \sigma \wedge \tau_a, S^g_n > x)$$
$$\leq (1 + \delta_a(x)) \sum_{n \geq 1} \mathbf{P}(n \leq \sigma \wedge \tau_a) \bar{F}(x + g(n)).$$

Again, given $\sigma \in \mathcal{T}$, $g \in \mathcal{G}_a$, and any a.s. finite stopping time $\varphi$, to prove (26), we may assume without loss of generality that $\varphi = m$ for some constant $m$. Since $\tau^m_a - m = \tau'_a$ where $\tau'_a = \min\{n \geq 1 : S'_n < an\}$, application of the result (27) to the random walk $\{S'_n\}_{n \geq 0}$ again given by $S'_n = S_{m+n} - S_m$, the stopping time $\sigma' = 0 \vee (\sigma - m)$ (for $\{S'_n\}$) and the function $g' \in \mathcal{G}_a$ now given by $g'(n) = g(m + n) - am$, gives, for $x > 0$,

$$\mathbf{P}(\exists n : 0 < n \leq (\sigma \wedge \tau^m_a) - m, S^g_{m+n} - S^{\bar{a}}_m > x)$$

$$= \mathbf{P}(\exists n : 0 < n \leq \sigma' \wedge \tau'_a, S'_n > x + g'(n))$$

(30)
$$\leq (1 + \delta_a(x)) \sum_{n \geq 1} \mathbf{P}(n \leq (\sigma \wedge \tau_a^m) - m) \bar{F}(x + g(m+n) - am)$$

$$\leq (1 + \delta_a(x)) \sum_{n \geq 1} \mathbf{P}(n \leq (\sigma \wedge \tau_a^m) - m) \bar{F}(x + g(m+n) - a(m+n)),$$

where the last line follows since $a > 0$. Now replace $n$ by $n - m$ in (30) to complete the proof as before. □

For any function $g$ on $\mathbb{Z}_+$ and any constant $a$, define the function $g^a$ by $g^a = g + \bar{a}$, so that, for each $n$, $g^a(n) = g(n) + an$.

We require also the following technical lemma.

LEMMA 2. *For any $\sigma \in \mathcal{T}$ and $g \in \mathcal{G}_0$, for all $0 < b < c$, and for all $x \geq 0$,*

$$H_\sigma^{g^b}(x) \geq H_\sigma^{g^c}(x) \geq \frac{b}{c} H_\sigma^{g^b}(x + c).$$

PROOF. The first inequality follows from the monotonicity of $\bar{F}$. To prove the second, for any $y \in \mathbb{R}_+$ define $\lceil y \rceil$ to be the least integer greater than or equal to $y$. Then, for $0 < b < c$ and all $y$,

$$c \lceil y \rceil \leq c(1 + y) \leq c + b \left\lceil \frac{c}{b} y \right\rceil,$$

and so

$$H_\sigma^{g^c}(x) = \int_0^\infty \mathbf{P}(\sigma \geq \lceil y \rceil) \bar{F}(x + c \lceil y \rceil + g(\lceil y \rceil)) \, dy$$

$$\geq \int_0^\infty \mathbf{P}\left(\sigma \geq \left\lceil \frac{c}{b} y \right\rceil\right) \bar{F}\left(x + c + b \left\lceil \frac{c}{b} y \right\rceil + g\left(\left\lceil \frac{c}{b} y \right\rceil\right)\right) dy$$

$$= \frac{b}{c} \int_0^\infty \mathbf{P}(\sigma \geq \lceil z \rceil) \bar{F}(x + c + b \lceil z \rceil + g(\lceil z \rceil)) \, dz$$

$$= \frac{b}{c} H_\sigma^{g^b}(x + c). \qquad \Box$$

PROOF OF THEOREM 2. We prove first (i). Fix $a > 0$ and define the sequence of a.s. finite stopping times $0 \equiv \rho^0 < \rho^1 < \rho^2 < \cdots$ for the process $\{S_n\}$ by, for $k \geq 1$,

$$\rho^k \equiv \rho_a^{\rho^{k-1}} = \rho^{k-1} + \min\{n \geq 1 : S_{\rho^{k-1} + n} - S_{\rho^{k-1}} > -an\}.$$

Note that $S_{\rho^k} > -a\rho^k$, $k \geq 0$, that is, that

(31) $$S_{\rho^k}^{-\bar{a}} > 0, \qquad k \geq 0.$$

For any $\sigma \in \mathcal{T}$, $g \in \mathcal{G}_0$, and for any $x > 0$, define the stopping time $\sigma_x$ by

$$\sigma_x = \sigma \wedge \min\{n : S_n^g > x\}.$$

Then

$$\mathbf{P}(M_\sigma^g > x) = \mathbf{P}\left(\bigcup_{k \geq 0} \{M_{\rho^k}^g \leq x; \exists n : \rho^k < n \leq \rho^{k+1}, \sigma \geq n, S_n^g > x\}\right)$$

$$= \sum_{k \geq 0} \mathbf{P}(M_{\rho^k}^g \leq x; \exists n : \rho^k < n \leq \rho^{k+1}, \sigma \geq n, S_n^g > x)$$

$$= \sum_{k \geq 0} \mathbf{P}(\exists n : \rho^k < n \leq \rho^{k+1}, \sigma_x \geq n, S_n^g > x)$$

$$\geq \sum_{k \geq 0} \mathbf{P}(\exists n : \rho^k < n \leq \rho^{k+1}, \sigma_x \geq n, S_n^g - S_{\rho^k}^{-\bar{a}} > x)$$

$$\geq (1 - \gamma_a(x)) \sum_{k \geq 0} \sum_{n \geq 1} \mathbf{P}(\rho^k < n \leq \rho^{k+1}, \sigma_x \geq n) \bar{F}(x + g(n) + an)$$

$$= (1 - \gamma_a(x)) \sum_{n \geq 1} \mathbf{P}(\sigma_x \geq n) \bar{F}(x + g(n) + an)$$

$$\geq (1 - \gamma_a(x)) \sum_{n \geq 1} (\mathbf{P}(\sigma \geq n) - \mathbf{P}(M_\sigma^g > x)) \bar{F}(x + g(n) + an),$$

where the fourth line in the above display follows by (31), while the fifth follows from Corollary 1(i) (with $\gamma_a$ as defined there). Since also $\sum_{n \geq 1} \bar{F}(x + g(n) + an) \leq \sum_{n \geq 1} \bar{F}(x + cn)$, it follows that

$$\mathbf{P}(M_\sigma^g > x)\left(1 + \sum_{n \geq 1} \bar{F}(x + cn)\right)$$

$$\geq (1 - \gamma_a(x)) \sum_{n \geq 1} \mathbf{P}(\sigma \geq n) \bar{F}(x + g(n) + an)$$

(32)
$$= (1 - \gamma_a(x)) H_\sigma^{g^a}(x)$$

$$\geq (1 - \gamma_a(x)) \frac{c}{c+a} H_\sigma^g(x + c + a),$$

where the last line above follows since the condition $g \in \mathcal{G}_c$ means that we can apply Lemma 2 to the function $g^{-c} \in \mathcal{G}_0$.

Observe that, as remarked in the Introduction, since the function $F$ is long-tailed, the function $H_\sigma^g$ is similarly long-tailed, with uniform convergence in the definition (1) over all $\sigma \in \mathcal{T}$ and $g \in \mathcal{G}_c$. Since also $\gamma_a(x) \to 0$ and $\sum_{n \geq 1} \bar{F}(x + cn) \to 0$, both as $x \to \infty$, it now follows from (32) that

$$\mathbf{P}(M_\sigma^g > x) \geq (1 - \gamma_a'(x)) H_\sigma^g(x)$$

for some positive function $\gamma'_a$, again independent of $\sigma$ and $g$, such that $\lim_{a \to 0} \lim_{x \to \infty} \gamma'_a(x) = 0$. The required lower bound (3) now follows, with uniformity over all $\sigma \in \mathcal{T}$ and $g \in \mathcal{G}_c$, as in the conclusion of the proof of part (i) of Lemma 1.

We now prove (ii). From the result (i), it is sufficient to show that

(33) $$\mathbf{P}(M^g_\sigma > x) \leq (1 + o(1)) H^g_\sigma(x) \qquad \text{as } x \to \infty,$$

uniformly over all stopping times $\sigma$ and all $g \in \mathcal{G}_c$. The proof is similar to, but simpler than, that of (i)—in particular, there is no need to define the stopping time $\sigma_x$. Fix $a \in (0, c)$ and define the sequence of a.s. finite stopping times $0 \equiv \tau^0 < \tau^1 < \tau^2 < \cdots$ for the process $\{S_n\}$ by, for $k \geq 1$,

$$\tau^k \equiv \tau_a^{\tau^{k-1}} = \tau^{k-1} + \min\{n \geq 1 : S_{\tau^{k-1}+n} - S_{\tau^{k-1}} \leq an\}.$$

Note that $S_{\tau^k} \leq a\tau^k$, $k \geq 0$, that is, that

(34) $$S^{\bar{a}}_{\tau^k} \leq 0, \qquad k \geq 0.$$

Then, for any stopping time $\sigma$, function $g \in \mathcal{G}_c$, and any $x > 0$,

$$\mathbf{P}(M^g_\sigma > x) \leq \sum_{k \geq 0} \mathbf{P}(\exists n : \tau^k < n \leq \tau^{k+1}, \sigma \geq n, S^g_n > x)$$

$$\leq \sum_{k \geq 0} \mathbf{P}(\exists n : \tau^k < n \leq \tau^{k+1}, \sigma \geq n, S^g_n - S^{\bar{a}}_{\tau^k} > x)$$

$$\leq (1 + \delta_a(x)) \sum_{k \geq 0} \sum_{n \geq 1} \mathbf{P}(\tau^k < n \leq \tau^{k+1}, \sigma \geq n) \bar{F}(x + g(n) - an)$$

$$= (1 + \delta_a(x)) \sum_{n \geq 1} \mathbf{P}(\sigma \geq n) \bar{F}(x + g(n) - an)$$

$$= (1 + \delta_a(x)) H^{g^{-a}}_\sigma(x),$$

where the function $\delta_a$ is as defined in Corollary 1(ii) above. Here the second line in the above display follows by (34), while the third follows from Corollary 1(ii). Hence, since again $g^{-c} \in \mathcal{G}_0$, it follows from Lemma 2 that, for $x \geq c$,

(35) $$\mathbf{P}(M^g_\sigma > x) \leq (1 + \delta_a(x)) \frac{c}{c - a} H^g_\sigma(x - c).$$

Note again that $H^g_\sigma$ is long-tailed, with uniform convergence in the definition (1) over all $\sigma \in \mathcal{T}$ and $g \in \mathcal{G}_c$. Hence, again arguing as in the conclusion of the proof of part (i), we obtain the required upper bound (33) with the required uniformity. $\square$

REMARK 3. The proof of Theorem 2 is close in spirit to that of Theorem 1 of [23].

**3. Comments, examples and counterexamples.** We give a number of examples and counterexamples, together with some commentary on the case where $\mathbf{P}(\sigma = \infty) > 0$. We continue to assume throughout that $F \in \mathcal{L}$ and that $F$ has mean zero.

In Examples 1–3, we show the importance of conditions on the functions $g$.

EXAMPLE 1. Here we show that, even for bounded stopping times, the functions $g$ cannot decrease too rapidly if we are to obtain uniform convergence over all $g$ in the conclusion (4). Suppose that $F \in \mathcal{S}$, and consider the stopping time $\sigma \equiv 2$. Consider also a sequence of functions $\{g_m\}_{m \geq 0}$ such that $g_m(1) = m$ and $g_m(2) = 0$ for all $m$. Then

$$\mathbf{P}(M_2^{g_m} > x) \geq \mathbf{P}(S_2^{g_m} > x) = 2(1 + o(1))\bar{F}(x) \qquad \text{as } x \to \infty,$$

while

$$H_2^{g_m}(x) = \bar{F}(x + m) + \bar{F}(x).$$

Hence, as in the discussion following Theorem 2, we obtain the conclusion (4), with $g = g_m$ for each fixed $m$. However, for any $\varepsilon > 0$ and for all sufficiently large $x$,

$$\liminf_{m \to \infty} \frac{\mathbf{P}(M_2^{g_m} > x)}{H_2^{g_m}(x)} \geq 2 - \varepsilon,$$

so that here the conclusion (4) does not hold with uniformity over all $m$.

EXAMPLE 2. Note that Theorems 1 and 2 extend to cover also functions $g$ which may take infinite values, provided that the definition (7) of $\mathcal{G}_c$ is interpreted as requiring that if, for any $n$, $g(n) = \infty$, then $g(n') = \infty$ for all $n' > n$. [A formal proof is given by replacing the stopping time $\sigma$ by $\sigma \wedge n$, where $n = \max\{n' : g(n') < \infty\}$ and using the existing results.]

In a continuation of the spirit of Example 1, suppose again that $F \in \mathcal{S}$ and consider now instead a function $g$ satisfying $g(1) = \infty$ and $g(2) = 0$. Fix $a > 0$ and define the stopping time $\sigma$ by $\sigma = 1$ if $\xi_1 \leq a$ and $\sigma = 2$ if $\xi_1 > a$. Then, as $x \to \infty$,

$$\mathbf{P}(M_\sigma^g > x) = \mathbf{P}(\xi_1 > a, \xi_1 + \xi_2 > x)$$
$$= \mathbf{P}(\xi_1 + \xi_2 > x) - \mathbf{P}(\xi_1 \leq a, \xi_1 + \xi_2 > x)$$
$$= (1 + o(1))(2\bar{F}(x) - \bar{F}(a)\bar{F}(x))$$
$$= (1 + o(1))(1 + \bar{F}(a))\bar{F}(x),$$

where the third line in the above display follows from the definition of subexponentiality and since also $F \in \mathcal{L}$. However,

$$H_\sigma^g(x) = \bar{F}(a)\bar{F}(x),$$

so that Theorem 1 will not extend to cover this case.

Now consider an alternative stopping time $\sigma'$ which is *independent* of $\{\xi_n\}_{n\geq 1}$ and has the same distribution as $\sigma$, that is, $\mathbf{P}(\sigma' = 1) = F(a)$ and $\mathbf{P}(\sigma' = 2) = \bar{F}(a)$. Then, as $x \to \infty$,

$$\mathbf{P}(M^g_{\sigma'} > x) = \bar{F}(a)\mathbf{P}(\xi_1 + \xi_2 > x) = (2 + o(1))\bar{F}(a)\bar{F}(x).$$

Since $H^g_{\sigma'}(x) = H^g_\sigma(x) = \bar{F}(a)\bar{F}(x)$, Theorem 1 again fails to extend to this case. However, this example also shows that, for this function $g$, the asymptotic distribution of the tail of $M^g_\sigma$ depends on $\sigma$ not just through its marginal distribution (as in the results of Theorems 1 and 2), but through the joint distribution of $\sigma$ and $\{\xi_n\}_{n\geq 1}$. See also [6] who consider a general function $g$ and a.s. constant stopping times.

EXAMPLE 3. In this example we show that, for a stopping time with unbounded support, and a function $g$ which increases too slowly, the tail of $\mathbf{P}(M^g_\sigma > x)$ may be heavier than that of $H^g_\sigma(x)$. Suppose that $g \equiv 0$ and that $\sigma$ is a random variable, independent of $\{S_n\}_{n\geq 0}$, such that $\mathbf{P}(\sigma > n) = (1+o(1))n^{-\alpha}$ as $n \to \infty$, for some $\alpha > 1$. Suppose also that the distribution $F$ has unit variance. Then

$$\mathbf{P}(M^g_\sigma > n) \geq \mathbf{P}(\sigma > n^2)\mathbf{P}(S_{n^2} > n) = (1 + o(1))cn^{-2\alpha} \qquad \text{as } n \to \infty,$$

where $c = \frac{1}{\sqrt{2\pi}}\int_1^\infty \exp\{-t^2/2\}\,dt$. We also have $H^g_\sigma(x) = \mathbf{E}\sigma \bar{F}(x)$ for all $x \geq 0$. Thus, if $F$ is additionally such that $\bar{F}(x) = o(x^{-2\alpha})$ as $x \to \infty$, then

$$\frac{\mathbf{P}(M^g_\sigma > x)}{H^g_\sigma(x)} \to \infty \qquad \text{as } x \to \infty.$$

The informal explanation here is that, for $g \equiv 0$, even moderate deviations contribute to the tail of $M^g_\sigma$. For more details on the asymptotics of $\mathbf{P}(M_n > x)$ as $n, x \to \infty$, see [5].

We now consider an example where the conditions of our main Theorem 2 do hold, and in which $\sigma < \infty$ a.s., but $\mathbf{E}\sigma = \infty$. In this case, when $F \in \mathcal{S}^*$ and $g \in \mathcal{G}_c$ for some $c > 0$, it follows, as in the derivation of (2), that $\bar{F}(x) = o(\mathbf{P}(M^g_\sigma > x))$ as $x \to \infty$, while, from Theorem 2(ii), we may deduce that $\mathbf{P}(M^g_\sigma > x) = o(\bar{F}^s(x))$. The example below shows that $\mathbf{P}(M^g_\sigma > x)$ may be of any order between $\bar{F}$ and $\bar{F}^s$.

EXAMPLE 4. Suppose that $F$ (which, as always, is assumed to have mean 0) is such that

$$\bar{F}(x) = (1 + o(1))K_2 x^{-\beta} \qquad \text{as } x \to \infty,$$

for some $K_2 > 0$ and $\beta > 1$. Then $F \in \mathcal{S}^*$ and

$$\bar{F}^s(x) = (1 + o(1))(\beta - 1)^{-1}K_2 x^{-\beta+1} \qquad \text{as } x \to \infty.$$

Consider any stopping time $\sigma$ with a tail distribution given by

$$\mathbf{P}(\sigma \geq n) = (1 + o(1))K_1 n^{-\alpha} \qquad \text{as } n \to \infty, \tag{36}$$

for some $K_1 > 0$ and $0 < \alpha < 1$. (E.g., since $F$ has finite variance, the distribution of the stopping time $\sigma = \min\{n : S_n > 0\}$ satisfies $\mathbf{P}(\sigma \geq n) = (1 + o(1))Kn^{-1/2}$ for some $K \in (0, \infty)$—see [13], Chapter 12.) Then $\mathbf{E}\sigma = \infty$ and, by Theorem 2(ii), for any $c > 0$ and as $x \to \infty$,

$$\mathbf{P}(M_\sigma^{\bar{c}} > x) = (1 + o(1)) \sum_{n \geq 1} \mathbf{P}(\sigma \geq n) \bar{F}(x + cn) \tag{37}$$

$$= (1 + o(1)) K_1 K_2 \sum_{n \geq 1} n^{-\alpha}(x + cn)^{-\beta} \tag{38}$$

$$= (1 + o(1)) K_1 K_2 \int_0^\infty t^{-\alpha}(x + ct)^{-\beta} dt$$

$$= (1 + o(1)) C x^{1-\alpha-\beta},$$

where

$$C = K_1 K_2 c^{\alpha-1} \int_0^\infty u^{-\alpha}(1 + u)^{-\beta} du,$$

and where (38) follows from (36) and (37) since the condition $\mathbf{E}\sigma = \infty$ implies that the contributions, as $x \to \infty$, of any finite number of the summands in (37) and (38) may be neglected.

In the case where $F$ has a Weibull distribution, that is, $\bar{F}(x) = (1 + o(1)) \exp(-x^\beta)$ as $x \to \infty$, for some $\beta \in (0, 1)$, then $\bar{F}^s(x) = (1 + o(1)) K_1 x^{1-\beta} \times \exp(-x^\beta)$ as $x \to \infty$. For the stopping time $\sigma$ as above and for $c > 0$, it follows similarly that

$$\mathbf{P}(M_\sigma^{\bar{c}} > x) = (1 + o(1)) K_2 x^{(1-\alpha)(1-\beta)} \exp(-x^\beta) \qquad \text{as } x \to \infty,$$

for some $K_2 > 0$.

We now discuss briefly the extent to which it is necessary that $\sigma$ should be a stopping time for the sequence $\{\xi_n\}_{n \geq 0}$ in order for our main results to hold.

In Example 5 we indicate briefly why some such condition is necessary.

EXAMPLE 5. Let $a > 0$ and define $\sigma = \min\{n : S_n > a\} - 1$. Then, for any nonnegative function $g$, $\mathbf{P}(M_\sigma^g > x) = 0$ for all $x \geq a$.

Now suppose again that $a > 0$ and consider the alternative stopping time $\sigma = \min\{n : \xi_n > a\} - 1$. Then by conditioning on each possible value of $\sigma$ and evaluating $\mathbf{P}(M_n^g > x | \sigma = n)$, one can straightforwardly show that $\mathbf{P}(M_\sigma^g > x) \leq c \exp(-\lambda x)$, for some constants $c > 0$ and $\lambda > 0$, so that here the distribution of $M_\sigma^g$ is again light-tailed, in contrast to the long tail of $H_\sigma^g$.

We now give, with some explanation, an example in which, although $\sigma$ is not a stopping time for the sequence $\{\xi_n\}_{n\geq 0}$, the equivalence (3) nevertheless holds.

EXAMPLE 6. Let $\{s_n\}_{n\geq 1}$ and $\{t_n\}_{n\geq 1}$ be two independent sequences of independent identically distributed random variables. Suppose that $s_1 \in \mathcal{S}$ with $\mathbf{E}s_1 = a$, and that $t_1 \geq 0$ a.s., with $\mathbf{E}t_1 = b > 0$. Let $T > 0$ be fixed, let $\eta = \min\{n : t_1 + \cdots + t_n > T\}$, and let $\sigma = \eta - 1$.

Let $c = b - a$ and define the sequence of independent identically distributed random variables $\{\xi_n\}_{n\geq 1}$, with distribution $F$, by $\xi_n = s_n - t_n + c$. Then, since $t_1$ is nonnegative and independent of $s_1 \in \mathcal{S}$, it follows easily that $\xi_1$ is tail-equivalent to $s_1$, and so also $\xi_1 \in \mathcal{S}$ and $\mathbf{E}\xi_1 = 0$. As usual, let $S_0 = 0$, $S_n = \sum_{i=1}^n \xi_i$, $n \geq 1$, be the random walk generated by the sequence $\{\xi_n\}_{n\geq 1}$. Then $M_\sigma^{\bar{c}} = \max_{0\leq n\leq \sigma} \sum_{i=1}^n (s_i - t_i)$ might, for example, be interpreted as the maximum loss to time $T$ of an insurance company with income at unit rate and a claim of size $s_n$ at each time $t_n$. Note that, clearly, $\mathbf{E}\exp(\lambda\sigma) < \infty$ for some $\lambda > 0$. Also $\sigma$ is not a stopping time for the random walk $\{S_n\}_{n\geq 0}$. However,

$$(39) \qquad \sup_{n\leq \sigma} \sum_{i=1}^n s_i - T \leq M_\sigma^{\bar{c}} \leq \sup_{n\leq \sigma} \sum_{i=1}^n s_i.$$

Since $T$ is fixed, $\sigma$ is independent of the sequence $\{s_n\}_{n\geq 1}$, and $s_1$ and $\xi_1$ are tail-equivalent, it follows from (39) and Theorem A 3.20 of [11], that, for any $c$,

$$(40) \quad \mathbf{P}(M_\sigma^{\bar{c}} > x) = (1 + o(1))H_\sigma^{\bar{c}}(x) = (1 + o(1))\mathbf{E}\sigma \bar{F}(x) \qquad \text{as } x \to \infty,$$

which is the equivalence (4) in this case. In the case $F \in \mathcal{S}^*$ and $c > 0$, we may go further and use Theorem 2(ii) of the present paper to obtain uniformity over all $T$ in the first equality in (40). See [18] for some further particular results on this model.

Note that the result follows here since $\sigma$ is a stopping time with respect to the sequence $\{s_n\}_{n\geq 1}$. In an intuitive sense (which might be made rigorous) the result also follows since, for each $n$, the event $\{\sigma \leq n\}$ is independent of the *tails* of the sequence $\xi_{n+1}, \xi_{n+2}\ldots$, and this is what is really required for our present results to hold.

Note also that the independence of the sequences $\{s_n\}_{n\geq 1}$ and $\{t_n\}_{n\geq 1}$ is vital. Consider instead a sequence $\{\xi_n\}_{n\geq 1}$ of independent identically distributed random variables with distribution $F \in \mathcal{S}^*$ and mean 0, and define the sequences $\{s_n\}_{n\geq 1}$ and $\{t_n\}_{n\geq 1}$ by $s_n = \max\{\xi_n, 0\}$ and $t_n = -\min\{\xi_n, 0\}$. Define $T$, $\eta$ and $\sigma$ as above. Then $\xi_\eta \leq 0$ a.s. and, for the random walk $\{S_n\}_{n\geq 0}$ generated by $\{\xi_n\}_{n\geq 1}$ and any $c > 0$, we have $M_\sigma^{\bar{c}} \equiv M_\eta^{\bar{c}}$.

Since $\eta$ is a stopping time for $\{S_n\}_{n\geq 0}$, it now follows from Theorem 2(ii) that

$$\mathbf{P}(M_\sigma^g > x) = (1+o(1))\mathbf{E}\eta \bar{F}(x) = (1+o(1))(\mathbf{E}\sigma + 1)\bar{F}(x) \qquad \text{as } x \to \infty.$$

EXAMPLE 7. Finally, we consider further the case of a stopping time $\sigma$ such that $p = \mathbf{P}(\sigma = \infty) > 0$. Recall that if $F \in \mathcal{S}^*$, then both $F \in \mathcal{S}$ and $F^s \in \mathcal{S}$. Provided only that $F^s \in \mathcal{S}$ (we do not here require our usual minimal assumption that $F \in \mathcal{L}$), and $m_F = 0$ as usual, then relatively straightforward arguments can be used to show that, in this case and for $c > 0$, the equivalence (4) continues to hold, and that, as $x \to \infty$,

$$\begin{aligned}
\mathbf{P}(M_\sigma^{\bar{c}} > x) &= (1+o(1))\mathbf{P}(\sigma = \infty)\mathbf{P}(M_\infty^{\bar{c}} > x) \\
&= (1+o(1))H_\sigma^{\bar{c}}(x) \\
&= \frac{(1+o(1))p}{c}\bar{F}^s(x).
\end{aligned} \tag{41}$$

However, under this weaker condition, we cannot expect any uniformity in either $\sigma$ or $c$.

In the case where $p = 1$ (i.e., $\sigma = \infty$ a.s.), the result (41) is the well-known theorem of Veraverbeke [19] referred to in the Introduction.

**Appendix** Recall that, for any stopping time $\sigma$ and nonnegative function $g$, the function $H_\sigma^g$ is defined by

$$H_\sigma^g(x) = \sum_{n\geq 1} \mathbf{P}(\sigma \geq n)\bar{F}(x+g(n)).$$

It is convenient to have a condition under which, for some purposes, we may replace the above sum by an integral.

Assume that, for $g \in \mathcal{G}_0$, the definition of the function $g$ is extended to all of $\mathbb{R}_+$ in such a way that $g$ continues to be increasing. For any such $g$, define the function $v^g$ on $\mathbb{R}_+$ by

$$v^g(x) = \sup_{n\geq 1} \frac{\bar{F}(x+g(n-1))}{\bar{F}(x+g(n))},$$

where $g(0) = 0$. For any stopping time $\sigma$ and $g \in \mathcal{G}_0$, define also the function $\widehat{H}_\sigma^g$ by

$$\widehat{H}_\sigma^g(x) = \int_0^\infty \mathbf{P}(\sigma > t)\bar{F}(x+g(t))\,dt.$$

Then, since $g$ is increasing and $\sigma$ is integer-valued, for all $x \in \mathbb{R}_+$,

$$\begin{aligned}
H_\sigma^g(x) &\leq \widehat{H}_\sigma^g(x) \\
&\leq \sum_{n\geq 1} \mathbf{P}(\sigma \geq n)\bar{F}(x+g(n-1)) \\
&\leq v^g(x)H_\sigma^g(x).
\end{aligned}$$

It follows, in particular, that if

(42) $$v^g(x) \to 1 \qquad \text{as } x \to \infty,$$

then also $H_\sigma^g(x) = (1 + o(1))\widehat{H}_\sigma^g(x)$ as $x \to \infty$.

Since $F \in \mathcal{L}$, the condition (42) holds for $g = \bar{c}$ [i.e., $g(n) = cn$] for any constant $c \geq 0$ (although observe that it does *not* hold with uniformity over all $c \geq 0$). More generally, the condition (42) holds for $g \in \mathcal{G}_0$ if $g(n) - g(n-1) \leq h(g(n))$ for some function $h$ satisfying (11).

**Acknowledgment.** The authors are most grateful to the referees for their very thorough reading of the paper and suggested improvements.

S. Foss
S. Zachary
Actuarial Mathematics and Statistics
Heriot–Watt University
Edinburgh EH14 4AS
Scotland
e-mail: s.zachary@ma.hw.ac.uk
e-mail: s.foss@ma.hw.ac.uk

Z. Palmowski
Mathematical Institute
University of Wrocław
pl. Grunwaldzki 2/4
50-384 Wrocław
Poland
and
Utrecht University
P.O. Box 80.010
3500 TA
Utrecht
The Netherlands
e-mail: zpalma@math.uni.wroc.pl